\documentclass{article}
\usepackage[english]{babel}
\usepackage{amsmath,amsfonts,amssymb,amscd,graphicx,color,amsthm,enumerate,calc} 
\usepackage{hyperref}

\def\OO{{\cal O}}
\def\PP{{\cal P}}

\newcommand{\Bc}{\ensuremath{{B^c}}}

\newcommand{\cobeta}{\ensuremath{{{\mathrm co}\beta}}}

\newcommand{\C}{{\ensuremath{\mathbb  C}}}
\newcommand{\Cbar}{\ensuremath{\overline{\mathbb  C}}}

\newcommand{\Cstar}{\ensuremath{{{\mathbb  C}^*}}}
\renewcommand{\d}{\ensuremath{\operatorname{d}}}

\newcommand{\diam}{\ensuremath{{\operatorname{diam}}}}

\newcommand{\Dbar}{\ensuremath{{\overline{\mathbb  D}}}}
\newcommand{\Dc}[1]{\ensuremath{{D^c_{#1}}}}
\newcommand{\Dcb}[1]{\ensuremath{{D^{c_b}_{#1}}}}
\newcommand{\Dcbar}[1]{\ensuremath{{\overline{D^c_{#1}}}}}
\newcommand{\Dcbbar}[1]{\ensuremath{{\overline{D^{c_b}_{#1}}}}}
\newcommand{\Dec}{{\ensuremath{\Delta}}}

\newcommand{\Decp}{{\ensuremath{\Delta'}}}
\newcommand{\DMp}[1]{\ensuremath{{D^\Mp_{#1}}}}
\newcommand{\DMpbar}[1]{\ensuremath{{\overline{D^\Mp_{#1}}}}}

\newcommand{\e}{\ensuremath{{\operatorname{e}}}}
\newcommand{\eps}{{\epsilon}}

\newcommand{\Ec}{\ensuremath{{E^c}}}

\newcommand{\Gc}{\ensuremath{{G^c}}}
\newcommand{\Gcb}{\ensuremath{{G^{c_b}}}}

\renewcommand{\H}{\ensuremath{\mathbb  H}}

\newcommand{\Hplus}{\ensuremath{{{\mathbb  H}_+}}}
\newcommand{\Hpq}{\ensuremath{{{H_{p/q}}}}}

\newcommand{\IMp}{\ensuremath{{I^{\Mp}}}}

\newcommand{\Log}{\ensuremath{\operatorname{Log}}}
\newcommand{\Lpq}{\ensuremath{L_{p/q}}}

\newcommand{\mapfromto}[3]{\hbox{\ensuremath{#1 : #2 \longrightarrow #3}}}
\renewcommand{\mod}{\ensuremath{{\operatorname{mod}}}}

\newcommand{\Mbrot}{\ensuremath{{\bf M}}}
\newcommand{\Mod}{\ensuremath{{\operatorname{mod}}}}
\newcommand{\Mone}{\ensuremath{{\bf M_1}}}
\newcommand{\Mp}{\ensuremath{{{M'}}}}
\newcommand{\Mpq}{\ensuremath{{{M_{p/q}}}}}
\newcommand{\N}{\ensuremath{\mathbb  N}}

\newcommand{\Q}{\ensuremath{\mathbb  Q}}

\newcommand{\sm}{\ensuremath{{\smallsetminus}}}

\newcommand{\Sen}{\ensuremath{{{\mathbb  S}^1}}}

\newcommand{\T}{\ensuremath{\mathbb T}}

\newcommand{\Vcb}{\ensuremath{{V^{c_b}}}}
\newcommand{\Vcn}{\ensuremath{{V_n^c}}}
\newcommand{\Vco}{\ensuremath{{V_1^c}}}
\newcommand{\Vcz}{\ensuremath{{V_0^c}}}

\newcommand{\whgamma}{\ensuremath{\widehat\gamma}}

\newcommand{\whB}{\ensuremath{\widehat B}}
\newcommand{\whBc}{\ensuremath{\widehat B^c}}

\newcommand{\whBMp}{\ensuremath{\widehat B^\Mp}}

\newcommand{\whD}{{\ensuremath{\widehat{D}}}}
\newcommand{\whDec}{{\ensuremath{\widehat{\Delta}}}}

\newcommand{\whkappa}{\ensuremath{\widehat \kappa}}

\newcommand{\whPi}{\ensuremath{\widehat \Pi}}

\newcommand{\whR}{\ensuremath{\widehat R}}
\newcommand{\whtau}{\ensuremath{\widehat\tau}}
\newcommand{\whtheta}{\ensuremath{\widehat\theta}}
\newcommand{\whT}{\ensuremath{\widehat T}}

\newcommand{\whXi}{{\ensuremath{\widehat{\Xi}}}}

\newcommand{\wtN}{\ensuremath{{\widetilde N}}}
\newcommand{\wtR}{\ensuremath{{\widetilde R}}}

\newcommand{\wtS}{\ensuremath{{\widetilde S}}}
\newcommand{\wtVc}{\ensuremath{{\widetilde V^c}}}
\newcommand{\wtVcb}{\ensuremath{{\widetilde V^{c_b}}}}
\newcommand{\wtVco}{\ensuremath{{\widetilde V_1^c}}}
\newcommand{\wtVcn}{\ensuremath{{\widetilde V_n^c}}}

\newcommand{\Wc}{\ensuremath{{W^c}}}

\def\eps{\epsilon}

\def\H{\mbox{$\mathbb H$}}
\def\C{\mbox{$\mathbb C$}}

\def\Z{\mbox{$\mathbb Z$}}
\def\Q{\mbox{$\mathbb Q$}}
\def\N{\mbox{$\mathbb N$}}

\def\Bottcher{B{\"o}ttcher}

\newtheorem{newthm}{Theorem}
\newtheorem{theorem}{Theorem}
\newtheorem{lemma}[theorem]{Lemma}
\newtheorem{proposition}[theorem]{Proposition}
\newtheorem{corollary}[theorem]{Corollary}
\newtheorem{remark}[theorem]{Remark}
\newtheorem{defthm}[theorem]{Definition and Theorem}
\newcommand{\ALIGN}{\begin{align*}}
\newcommand{\ENDALIGN}{\end{align*}}
\newcommand{\ENUM}{\begin{enumerate}}
\newcommand{\ENUMa}{\begin{enumerate}[a.]}
\newcommand{\ENUMi}{\begin{enumerate}[i)]}
\newcommand{\ENDENUM}{\end{enumerate}}
\newcommand{\ITMZ}{\begin{itemize}}
\newcommand{\ENDITMZ}{\end{itemize}}
\newcommand{\REFEQN}[1] { \begin{equation}\label{#1} }
\newcommand{\ENDEQN}{\end{equation}}
\newcommand{\THM}{\begin{theorem}}
\newcommand{\REFEXA}[1] { \begin{example}\label{#1} }
\newcommand{\ENDEXA}{\end{example}}
\newcommand{\REM}{ \begin{remark}}
\newcommand{\ENDREM}{\end{remark}}
\newcommand{\REFTHM}[1] { \begin{theorem}\label{#1} }
\newcommand{\RREFTHM}[2] { \begin{theorem}[#1]\label{#2} }
\newcommand{\ENDTHM}{\end{theorem}}
\newcommand{\REFNTH}[1] { \begin{newthm}\label{#1} }
\newcommand{\ENDNTH}{\end{newthm}}
\newcommand{\REFPROP}[1]{\begin{proposition}\label{#1} }
\newcommand{\RREFPROP}[2]{\begin{proposition}[#1]\label{#2} }
\newcommand{\PROP}{\begin{proposition}}
\newcommand{\ENDPROP}{\end{proposition} }
\newcommand{\REFDEF}[1]{\begin{definition}\label{#1} }
\newcommand{\DEF}{\begin{definition}}
\newcommand{\ENDDEF}{\end{definition} }
\newcommand{\REFLEM}[1]{\begin{lemma}\label{#1} }
\newcommand{\RREFLEM}[2]{\begin{lemma}[#1]\label{#2} }
\newcommand{\LEM}{\begin{lemma}}
\newcommand{\ENDLEM}{\end{lemma} }
\newcommand{\REFCOR}[1]{\begin{corollary}\label{#1} }
\newcommand{\COR}{\begin{corollary}}
\newcommand{\ENDCOR}{\end{corollary}}
\newcommand{\RMRK}{\begin{remark}}
\newcommand{\ENDRMRK}{\end{remark}}
\newcommand{\REFDEFTHM}[1] { \begin{defthm}\label{#1} }
\newcommand{\ENDDEFTHM}{\end{defthm}}

\newcommand{\corref}[1]{Corollary~\ref{#1}}

\newcommand{\figref}[1]{Fig.~\ref{#1}}
\newcommand{\lemref}[1]{Lemma~\ref{#1}}
\newcommand{\thmref}[1]{Theorem~\ref{#1}}
\newcommand{\propref}[1]{Proposition~\ref{#1}}

\newcommand{\PROOF}{\begin{proof}}
\newcommand{\ENDPROOF}{\end{proof}}

\title{Carrots for dessert}
\author{Carsten Lunde Petersen and Pascale Roesch}

\begin{document}

\maketitle

\begin{abstract}
Carrots for dessert is the title of a section of the paper `On polynomial-like mappings', 
\cite{polylikemaps}. In that section Douady and Hubbard define a notion of 
dyadic carrot fields of the Mandelbrot set and more generally for Mandelbrot like families 
(for a precise statement see below). 
They remark that such carrots are small when the dyadic denominator is large, but they 
do not even try to prove a precise such statement. 
In this paper we formulate and prove a precise statement of asymptotic 
shrinking of dyadic Carrot-fields around $\Mbrot$. The same proof carries readily over to 
show that the dyadic decorations of copies $\Mp$ of the Mandelbrot set $\Mbrot$ inside $\Mbrot$ 
and inside the parabolic Mandelbrot set $\Mone$ shrink to points when the denominator diverge to $\infty$.
\end{abstract}

\section*{Introduction}
For $c\in\C$ let $Q_c(z) = z^2 +c$ and let $J_c$ and $K_c$ denote respectively 
the Julia set and the filled Julia set for $Q_c$. Denote by $\Mbrot$ the Mandelbrot set
$$
\Mbrot = \{c\in\C\;|Q_c^n(0) \underset{n\to\infty}{\nrightarrow}\infty\}.
$$
Similarly for $B\in\C$ let $g_B(z) = z + 1/z + B$. 
Then each $g_B$ has a parabolic fixed point at $\infty$ with multiplier $1$ and 
$g_B$ is conjugate to $g_{-B}$ via $z\mapsto -z$. 
The parabolic Mandelbrot set $\Mone$ is the set 
\begin{equation}\label{Mone}
\Mone = \{A\in\C| \textrm{either } g_{\sqrt{A}}^n(-1)\underset{n\to\infty}{\nrightarrow}\infty 
\textrm{ or } g_{\sqrt{A}}^n(1)\underset{n\to\infty}{\nrightarrow}\infty\}.
\end{equation}
Let $T$ be a closed triangle in the right halfplane $\Hplus$ union $\{0\}$ bounded by lines through the origin 
and a non horizontal line in such a way that the $i2\pi\Z$ translates are disjoint. 
Let $\whDec_0$ be the image of $T$ under $z\mapsto\e^{z}$. 
Then by construction $\whDec_0$ is simply connected and $Q_0^{-1}(\whDec_0)$ has two connected 
components one, which is a subset of $\whDec_0$ and another one $\whDec_{1/2}$ containing 
$-1=\e^{i2\pi/2}$. 
Define recursively $\whDec_{p/2^n}$ as the connected component 
of $Q_0^{-n}(\whDec_0)$ containing $\exp(i2\pi p/2^n)$. 
The sets $\whDec_0$ and $\whDec_{p/2^n}$, $0<p<2^n$, $p$ odd and $n\in\N$ 
are disjoint. Together they form a ``dyadic Carrot field'' $\whDec$ of $\Dbar$: 
$$
\whDec=\whDec_0\cup \bigcup_{n\geq 1}\bigcup_{0<p<2^n} \whDec_{p/2^n}.
$$

The degenerate version of such a carrot field is a ``dyadic stick field'' defined 
similarly, but with $T = [1,t]$ for some $t>1$. 
We shall in the following denote by carrot field any possibly degenerate carrot field.

Let {\mapfromto {\Psi} {\Cbar\sm\Dbar} {\Cbar\sm\Mbrot}} denote the Douady-Hubbard 
uniformizing parameter. That is $\Psi$ is biholomorphic, tangent to the identity at $\infty$ and 
its inverse $\Phi$ is given by $\Phi(c) = \phi_c(c)$, where $\phi_c$ denotes 
the {\Bottcher}-coordinate of $Q_c$ at $\infty$. 
We shall use also the Green's functions for $\Mbrot$ and $K_c$, i.e. 
the subharmonic functions $g_\Mbrot(c) = \log^+|\phi_c(c)|$ and
$$
g_c(z) = \lim_{n\to\infty} \frac{1}{2^n}\log^+|Q_c^n(z)|.
$$ 
Moreover we shall use the notation $E_\Mbrot(h)$ and $E_c(h)$ for the equipotentials 
for $g_\Mbrot$ and $g_c$ of level $h\geq 0$. 
Similarly we shall use the notation $F_\Mbrot(h)$ and $F_c(h)$ for the closed filled equipotentials 
of level or height $h$:
$$
F_\Mbrot(h) = \{c|g_\Mbrot(c)\leq h\}, \qquad F_c(h) = \{z|g_c(z)\leq h\}.
$$
The external ray of argument $\theta$ for  $\Mbrot$ or $K_c$ 
is the field line of $g^\Mbrot$ or $g^c$, 
which is asymptotic to the halfline $\exp(t+i2\pi\theta)$ at $\infty$. 

By the Douady-Hubbard landing theorem for rational external rays of $\Mbrot$, 
$\Psi$ has a continuous extension along any external ray of rational argument. 
In particular along the rays $R_\theta$ with dyadic arguments $\theta=p/2^n$. 
By construction each connected component $\whDec_\theta$ of $\whDec$ is contained 
in some Stolz angle measured from its vertex $\e^{i2\pi\theta}\in\Sen$. 
It thus follows that the radial extension of $\Psi$ along $R_\theta$ defines 
a continuous extension of $\Psi$ to $\Psi(\whDec_\theta)$ for each dyadic $\theta$. 
Hence $\Psi(\whDec)$ is well defined. 
Write $\Dec=\Psi(\whDec)$ and $\Dec_{p/2^n} = \Psi(\whDec_{p/2^n})$. 

Then the carrot (or stick) decorated Mandelbrot set is $\Mbrot\cup\Dec$, 
where $\Dec=\Psi(\whDec)$ is any (possibly degenerate) carrot field.

We can also at least partially transport $\whDec$ to the dynamical plane of $Q_c$ and 
thus obtain $\Dec^c=\phi_c^{-1}(\whDec)$, where we use for $\phi_c^{-1}$ the maximal radial extension. 
We can then view the carrots $\Dec$ of $\Mbrot$ as the set of parameters for which 
$c$ belongs to the corresponding carrot $\Dec^c$ of the filled Julia set $K_c$. 

With this terminology the Theorem of shrinking of dyadic carrots of $\Mbrot$ is 
\REFTHM{shrinkingcarrots}
For any (possibly degenerate) dyadic carrot field $\Dec$ of $\Mbrot$
$$
\lim_{n\to\infty}\diam(\Dec_{p/2^n}) = 0.
$$
\ENDTHM
We shall refer to any of the sets $\Dec_{p/2^n}$ as a dyadic carrot of $\Mbrot$.

An easy adaptation of our proof shows that dyadic is not essential. That is if 
$\whDec$ is a carrot field, where instead $\whDec_0$ is any finite collection of disjoint (possibly degenerate) triangles attached to periodic orbits for $Q_0$ and $\whDec$ is obtained by iterated pull back as above. 
Then the corresponding version of \thmref{shrinkingcarrots} still holds. 

Let $M'$ with period $k$ denote a copy of $\Mbrot$ inside $\Mbrot$ 
or a copy of $\Mbrot$ inside the Parabolic Mandelbrot set $\Mone$. 
Let $\theta\pm$ be the arguments of the pair of external rays (parabolic external rays if $M'\subset\Mone)$ 
co-landing at the root $c'_0$ of the principal hyperbolic component $H'$ for $M'$. 
We let {\mapfromto {\chi_{\Mp}} {\Mp} \Mbrot} denote the Douady-Hubbard straightening map 
(for a definition see \cite[Chap. II, l-4]{polylikemaps}). 

Let $I = I(H') = I(\Mp) = [\theta_-,\theta_+]$ be the tuning interval for $\Mp$ (or equivalently for $H'$) 
and let $\whtheta_+<\whtheta_-\in I$ be the points such that 
each of the subintervals $I_0=[\theta_-,\whtheta_+]$ and $I_1=[\whtheta_-,\theta_+]$ map diffeomorphically 
onto $I$ under $\sigma^k$, where $\sigma(\theta) = 2\theta\;\Mod\; 1$. 
Let $\IMp$ denote the corresponding $\sigma^k$-invariant Cantor set and let 
{\mapfromto {\kappa=\kappa_\Mp} {\IMp} {\Sigma_2}} denote the conjugacy 
of {\mapfromto {\sigma^k} \IMp \IMp} to the shift on $\Sigma_2={\{0,1\}}^{\N}$ with 
$\kappa(\theta_-) = \overline{0}$ and $\kappa(\theta_+) = \overline{1}$. 
Then the pair of rays with arguments $\whtheta_\pm$ coland at the principal tip 
$c_{1/2}'=\chi_{M'}^{-1}(\Psi(1/2))$ of $\Mp$. 
The sector $W_{1/2}'$ bounded by these rays and disjoint from $\Mp$
is called the principal wake of $\Mp$ and the intersection $\Decp_{1/2} := \overline{W_{1/2}'}\cap \Mbrot$ 
is called the $1/2$ dyadic decoration of $\Mp$. 
More generally for $p$ odd with binary representation $p=\eps_1\ldots\eps_n$, $\eps_n\not=0$ the dyadic 
number $p/2^n$ has two binary representations $0.\eps_1\ldots\eps_n\overline{0}$ and $0.\eps_1\ldots\eps_{n-1}0\overline{1}$. 
According to the Douady tuning algorithm $\theta_{p/2^n}^-=\kappa^{-1}(\eps_1\ldots\eps_n\overline{0})$ and 
$\theta_{p/2^n}^+=\kappa^{-1}(\eps_1\ldots\eps_{n-1}0\overline{1})$ are the two endpoints of 
a complementary interval of $\IMp$. 
Moreover the corresponding external rays of $\Mbrot$ co-land at the relatively dyadic tip 
$c'_{p/2^n}=\chi_{M'}^{-1}(\Psi(p/2^n))$ of $\Mp$ and for any parameter $c\in\Mp$ the 
corresponding dynamical rays co-land on a point, which is preperiodic to the relative $\beta$ fixed point. 
The $p/2^n$-wake $W'_{p/2^n}$ and the dyadic decoration $\Decp_{p/2^n} := \overline{W'_{p/2^n}}\cap \Mbrot$ 
of $\Mp$ are defined similarly as the $1/2$ wake and decoration. 
Denote by $W_0'$ the sector bounded by the rays of arguments $\theta_\pm$ 
and not containing $\Mp$ and let $\Decp_0 = \overline{W_0'}\cap\Mbrot$. 
Note that for each $p/2^n$ the root $c'_{p/2^n}$ of the corresponding wake or limb 
is the only point of intersection between $\Mp$ and (the closure of the) wake or limb. 
Note also that any two wakes are disjoint. As above we write 
$$
\Decp=\Decp_0\cup \bigcup_{n\geq 1}\bigcup_{0<p<2^n} \Decp_{p/2^n}.
$$
(For $\Mp$ a copy of $\Mbrot$ inside $\Mone$ we use parabolic rays.).

Then 
\begin{theorem}[Douady-Hubbard, Yoccoz]\label{dyadicMdecomposition}
For any copy $\Mp$ of $\Mbrot$ in $\Mbrot$ :
$$
\Mbrot = \Mp\cup\bigcup_{n\geq 0}\bigcup_{p/2^n\in\Q} \Decp_{p/2^n}.
$$
\end{theorem}

\PROOF
The copy $\Mp$ of $\Mbrot$ belongs to the limb $L_{p'/q'}^\Mbrot$ of 
the central hyperbolic component $H_0$ of $\Mbrot$, for some 
$p',q'\in\N$ with $(p',q')=1$. 
Let $c\in\Mp$ be the center of $\Mp$, i.e.  $Q_c^k(c) = c$, where $k$ is the period of $\Mp$.
Let $P_n$, $n\in\N$ be the level $n$ ($p'/q'$)-Yoccoz puzzle piece containing the critical value $c$ 
of $Q_c$ and let $\PP_n$ denote ($p'/q'$)-Parameter Yoccoz puzzle piece containing the parameter $c$. 
Then for each $n$ the map $\Psi\circ\phi_c(z)$ restricted to $\partial P_n\sm J_c$ 
extends to a homeomorphism of $\partial P_n$ onto $\partial\PP_n$ preserving argument and 
potential. 
Also for each $n$ the closed puzzle piece $\overline{P}_n$ contains the ends from potential 
$2^{-n}$ and down of the external rays with arguments in $\IMp$. 
Hence the same holds for the corresponding parameter rays and $\overline{\PP}_n$. 
It follows that any other level $n$ parameter puzzle piece as well as 
$\Mbrot\sm L_{p'/q'}^\Mbrot$ is contained in one of the relatively dyadic wakes $\Decp_{p/2^m}$ of $\Mp$. 
The theorem then follows from Yoccoz parameter puzzle theorem for renormalizable 
parameters, which states that 
$$
\Mp = \bigcap_{n\geq 0} \overline{\PP}_n.
$$
\ENDPROOF

\begin{theorem}\label{dyadicMonedecomposition}
For any copy $\Mp$ of $\Mbrot$ in $\Mone$ :
$$
\Mone = \Mp\cup\bigcup_{n\geq 0}\bigcup_{p/2^n\in\Q} \Decp_{p/2^n}
$$
\end{theorem}
\PROOF
Completely analogous to the above.
\ENDPROOF

The Shrinking decorations Theorem for strict copies $\Mp$ of $\Mbrot$ in $\Mbrot$ or $\Mone$ 
can then be stated as
\REFTHM{shrinkingcarrotsp}
For any strict copy $\Mp$ of $\Mbrot$ in $\Mbrot$ or in $\Mone$
$$
\lim_{n\to\infty}\diam(\Decp_{p/2^n}) = 0.
$$
\ENDTHM
The two theorems \thmref{shrinkingcarrots} and \thmref{shrinkingcarrotsp} have 
very similar proofs, the proof of the first being slighly more complicated. 
We shall detail the proof of the first and sketch the difference to the proof 
of the second. 

Dzmitry Dudko presents a different and independent proof of the Shrinking decorations Theorem 
for strict copies $\Mp$ of $\Mbrot$ in $\Mbrot$ in \cite{Dudko}. His statement includes more generally 
strict copies of the Multibrot set inside the Multibrot set of the same degree. 
The proof we give here would also easily extend to the Multibrot case.

\section*{Proofs}
\subsubsection*{First reduction: Independence on $T$.}
Going back to the initial setting of possibly degenerate carrot fields decorating $\Mbrot$. 
We shall first show that the proof of \thmref{shrinkingcarrots} can be reduced to 
considering only one particular stick-field. 

Indeed let $T^1$ and $T^2$ be any two possibly degenerate triangles in $\Hplus\cup\{0\}$. 
and let $\check{T}^i=T^i\sm\{0\}$ for $i=1,2$. 
Then there exists $\delta>0$ such that $\check{T}^1$ is contained in a 
hyperbolic $\delta$-neighbourhood of $\check{T}^2$ in $\Hplus$ and vice versa. 
As {\mapfromto {\exp} \Hplus {\C\sm\Dbar}} and {\mapfromto {Q_0} {\C\sm\Dbar} {\C\sm\Dbar}} 
are hyperbolic isometries the same statement holds for $\whDec_\theta^i$, $i=1,2$ and 
$\theta=p/2^n$ any dyadic. 
By elementary estimates on hyberbolic metrics, there exists $k=k(\delta)>1$ such that 
for any univalent map {\mapfromto \psi {\C\sm\Dbar} \C}, tangent to the identity at infinity, 
and any dyadic $\theta$ 
$$
\frac{1}{k} \leq 
\frac{\diam(\psi(\whDec_\theta^1\sm\{\e^{i2\pi\theta}\}))}{\diam(\psi(\whDec_\theta^2\sm\{\e^{i2\pi\theta}\}))}
\leq k
$$
in particular 
$$
\frac{1}{k} \leq 
\frac{\diam(\Dec_\theta^1)}{\diam(\Dec_\theta^2)},
\leq k
$$
where $\diam(\cdot)$ denotes euclidean diameter. 

Hence to prove \thmref{shrinkingcarrots} it suffices to consider a particular stick field, 
say the field for $t=1/2$, which is what we shall do.

\subsubsection*{The toy, but key argument}
To set the scene let us however consider first the toy example, where we replace the interval 
$T = [0,1/2]$ defining $\whDec_0$ by a compact subset of $\H$, whose $i2\pi\Z$ translates are disjoint, 
i.e.~whose projection to $\C\sm\Dbar$ does 
not separate $\Dbar$ from $\infty$, say $\whDec_0=\exp([1/4,1/2])$. 
This completely trivialises the problem by considerations 
on the comparison of hyperbolic and euclidean distance similar to above:
In this simpler case the set  $\whDec_0$ and thus also $\Dec_0$ has finite hyperbolic diameter $diam$ 
and moreover this bound on the hyperbolic diameter of $\whDec_0$ is an upper bound on the hyperbolic diameter 
of any of the dyadic carrots $\Dec_{p/2^n}\subset\Dec$. 
Hence the euclidean diameter of any such dyadic carrot is bounded uniformly from above 
by a universal constant $k = k(diam)$ times the euclidean distance between 
$\Dec_{p/2^n}$ and $\Mbrot$. 
Since the later tends to zero uniformly as $n\to\infty$ we have in the toy case 
$$
\limsup_{n\to\infty} \diam(\Dec_{p/2^n}) = 0.
$$

We shall see that, this is effectively what happens, if we restrict our attention to 
any renormalization copy $\Mp$ of $\Mbrot$. 
However the decorations around $\Mp$ are not the image under a holomorphic map 
of a union of compact sets all of which are isometric copies of an inital connected component. 
Hence we need to device other means of making hyperbolic estimates. 
To this end we shall use the fact that if $K\subset V'\subset V\subset U$, 
with $U$ a hyperbolic domain and $\mod(V\sm\overline{V'})\geq \delta>0$, 
then the hyperbolic diameter of $K$ in $U$ satisfies $\diam_U(K) \leq d(\delta)$. 
And we shall use the observation by Shishikura, 
that holomorphic motions can be used to transfer bounds 
for (locally) persistent annuli in dynamical space to bounds for corresponding 
annuli in parameter space (see \cite{Roesch}). 

With this in mind let us proceed to the proof of \thmref{shrinkingcarrots}. 
Then as mentioned above \thmref{shrinkingcarrotsp} will follow by using the same proof.

\subsubsection*{Proof of \thmref{shrinkingcarrots}.}
We prove the following result : 
\begin{proposition}\label{shrinkingofconvgcarrots}
For any sequence ${\{\Dec_k = \Dec_{p_k/2^{n_k}}\}}_{k\in\N}$, $n_{k+1}>n_k$ of carrots for $\Mbrot$ 
with roots $c_k$:
$$
c_k\underset{k\to\infty}\longrightarrow c_\infty \Longrightarrow 
\diam(\Dec_k)\underset{k\to\infty}\longrightarrow 0.
$$
\end{proposition}

\begin{remark} \thmref{shrinkingcarrots} is an easy corollary of this proposition 
by compactness of the Mandelbrot set. The details are left to the reader.

\end{remark}

\subsubsection*{Setup for the proof of \propref{shrinkingofconvgcarrots}.}
Let  ${\{\Dec_k = \Dec_{p_k/2^{n_k}}\}}_{k\in\N}$, $n_{k+1}>n_k$ be an arbitrary but fixed sequence 
of carrots for $\Mbrot$ with roots $c_k$ converging to $c_\infty$. 
Then first of all $c_\infty\in\partial\Mbrot$.

\vskip 1em
We shall use the Levin-Yoccoz parameter space inequality and 
Yoccoz theorem on local connectivity of $\Mbrot$ at Yoccoz parameters, 
i.e. parameters $c$, for which $Q_c$ is not (infinitely) renormalizable 
and has all periodic points repelling. 
For the version of \propref{shrinkingofconvgcarrots} leading to a proof of 
\thmref{shrinkingcarrotsp} the simpler Yoccoz (rather than Levin-Yoccoz) 
parameter space inequality suffices, 
but for \propref{shrinkingofconvgcarrots} we need the extension due to Levin:

\THM[The Yoccoz-Levin Dynamical Inequality]
Let $H$ be any hyperbolic component of $\Mbrot$ of period $k$. 
Let $p/q$ be any non zero reduced rational and let $W_{p/q}^H$ denote 
the relative $p/q$ wake of $H$, bounded by parameter rays with arguments 
$0<\eta_-<\eta_+<1$. 
For any $c\in W_{p/q}^H$ let $\lambda$ denote the multiplier of the repelling 
$k$-periodic common landing point $\alpha'$ of the $kq$ periodic rays $R_{\eta_\pm}^c$. 
Then $\alpha'$ has combinatorial rotation number $p/q$ and 
 $\lambda$ has a logarithm $\Lambda$ such that:
$$
|\Lambda -p/q2\pi i| \leq \frac{2k\log2\cos\theta}{q}\frac{\pi}{\omega(c)}, 
$$
where $\theta\in\;]-\pi/2,\pi/2[$ is the argument of $\Lambda -p/q2\pi i$ 
and $\omega(c)$ is the angle of vision of the interval $i2\pi[\eta_-,\eta_+]$ 
from $\Log \phi_c(c)\in\{z=x+iy| 0< y < 2\pi\}$.
\ENDTHM
\newcounter{YoccozLevincounter}\setcounter{YoccozLevincounter}{\value{theorem}}
\PROOF
Levin proved the fixed point case $k=1$ in \cite[TH. 5.1]{Levin}, 
the general case is similar. For completeness we give a proof in the Appendix, 
page~\pageref{YoccozLevindynamicalproof}. 
\ENDPROOF

\PROP
Let $H$ be a period $k$ hyperbolic component of $\Mbrot$ with tuning interval $I(H) = [\theta_-,\theta_+]$ 
and let $p'/2^{m'}\in I(H)$, $1\leq m' <k$ be the dyadic with the smallest denominator. 
For any irreducible rational $p/q$, let  $0<\eta_-<\eta_+<1$ be 
the arguments of the co-landing parameter rays bounding $W_{p/q}^H$ 
and let $p/2^m\in [\eta_-,\eta_+]$ be the dyadic with the smallest denominator. 
We have 
$$
2^{-kq} \leq \eta_+-\eta_-\qquad\text{and}\qquad m = m'+k(q-2).
$$
\ENDPROP
\PROOF
As $\eta_-<\eta_+$ are periodic of exact period $kq$, we have 
$\eta_+-\eta_- \geq 1/(2^{kq}-1) > 2^{-kq}$. For the second inequality let $\Mp$ denote the copy of 
$\Mbrot$ with $H$ as central hyperbolic component. 
Let $\theta_-<\theta_+\in\IMp$ denote the arguments of the parameter rays colanding at the root of $\Mp$. 
Let $I=[\theta_-,\theta_+]\supset I_0, I_1$,  $\IMp\subset I_0\cup I_1$ 
and {\mapfromto \kappa \IMp {\Sigma_2}} be as above and write $\pi$ for the binary projection 
of $\Sigma_2$ onto $\T$ and set $\whkappa=\pi\circ\kappa$. 
Then $\tau_\pm = \whkappa(\eta_\pm)$ are the arguments of the parameter rays co-landing 
at the root of the wake $W_{p/q}^{H_0}$. 
It is well known that the intervals $\sigma^j([\tau_-,\tau_+])$, $0\leq j < q$ are interiorly disjoint and 
injective images. Moreover $0\in\sigma^{(q-1)}([\tau_-,\tau_+])$ and thus 
$1/2\in\sigma^{(q-2)}([\tau_-,\tau_+])$. 
Consequently $\sigma^{k(q-2)}$ maps $[\eta_-,\eta_+]$ injectively into $I$. 
Morever $I \supset \sigma^{k(q-2)}([\eta_-,\eta_+])\supset (I\sm (I_0\cup I_1))$. 
Let $p'/2^{m'}\in I$ be the dyadic with smallest denominator then 
$p'/2^{m'}\in I\sm (I_0\cup I_1)$ and $1\leq m'\leq k$. 
Thus $m = m'+k(q-2)$.
\ENDPROOF

\REFCOR{angleofdecorationview}
For any $c\in W_{p/q}^H\cap(\Mbrot\cup\Dec)$ 
the angle of vision $\omega(c)$ of $i2\pi[\eta_-,\eta_+]$ from $\Log(\phi_c(c))$ is bounded from below by 
$$
\arctan(2\pi2^{m'-2k}).
$$
\ENDCOR
\PROOF
The angle is bounded from below by the angle obtained, 
when $c$ belongs to one of the two bounding rays of $W_{p/q}^H$:
\begin{align*}
\arctan(2\pi(\eta_+-\eta_-)/\log|\phi_c(c)|) &\geq \arctan(2\pi 2^{-kq}/2^{-(m'+k(q-2))})\\
&= \arctan(2\pi2^{m'-2k}).
\end{align*}
\ENDPROOF

\THM[The Yoccoz-Levin Parameter Inequality]\label{yoccozinequality}
For any hyperbolic component $H$ of $\Mbrot$ there exists a constant $C=C_H>0$ 
such that for any relative $p/q$ wake $W_{p/q}^H$ 
$$
\diam(W_{p/q}^H\cap(M\cup\Dec)) \leq \frac{C}{q}.
$$
\ENDTHM
\PROOF
For $h=0$ i.e. for the limbs $\Mbrot\cap W_{p/q}^H$, 
this is essentially proved by Hubbard in \cite{Hubbard}, 
except that he confuses the direction of the square root from primitive hyperbolic components 
and obtains an inequality with $C/\sqrt{q}$ instead of $C/q$ in the primitive case. 
Whereas the bounds actually gives $C/q^2$ assymptotically when $p/q$ tend to $0$ or $1$. 
For the extension we use the Levin-Yoccoz dynamical inequality above instead of the Yoccoz inequality. 
By \corref{angleofdecorationview} the angle $\omega(c)$ for $c\in W_{p/q}^H\cap(M\cup\Dec)$ 
is bounded from below by the angle $\omega_H=\arctan(2\pi2^{m'-2k})$. 
The argument is then identical to the argument in Hubbards paper \cite{Hubbard}, 
except using the Levin-Yoccoz dynamical inequality with the fixed value $\omega_H$. 
Thus asymptotically for $q$ large we can take
$$
C_H = \frac{\pi}{\omega_H}C_H^{\textrm{Yoccoz}}
$$
where $C_H^{\textrm{Yoccoz}}$ is the corresponding assymptotic value for Yoccoz parameter inequality.
\ENDPROOF

Let $H_0$ denote the central hyperbolic component of $\Mbrot$ and for 
$p/q$ an irreducible rational let  $W_{p/q}^{H_0}$ denote the $p/q$ wake of $H_0$ and 
$\Lpq^{H_0}=W_{p/q}^{H_0}\cap\Mbrot$ the $p/q$ limb. 

\RREFTHM{Yoccoz}{yoccozpuzzle}
For any $p/q$ and any $c\in\Lpq^{H_0}$ there are two possibilities; 
either $c$ is not renormalizable and the $p/q$ parameter puzzle pieces containing $c$ 
nests down to $c$, or $c$ is at least once renormalizable, say with first renormalization period $k$  
and there is a first level $n$ such that for the dynamical puzzle pieces $P_n =: U$ and $P_{n+k} =: U'$,  
{\mapfromto {Q_c^k} U {U'}} is quadratic like with connected filled-in Julia set. (fattening $U$ and $U'$ if $k=q$.)
\ENDTHM
For a proof see \cite{Hubbard}.

\subsubsection*{The second reduction : reduction to renormalizable  $c_\infty$.}
Let us first apply the Yoccoz-Levin parameter inequality. 
This gives a constanct $C>0$ such that for all $p/q$ 
$$
\diam(W_{p/q}^{H_0}\cap(\Mbrot\cup \Dec)) \leq \frac{C}{q}.
$$
The sequence  ${\{\Dec_k \}}$ of carrots is included in a sequence of Wakes $W_{p_k/q_k}^{H_0}$. 
Hence it follows that if $q_k$ tends to $\infty$ the diameter of $\Dec_k$ tends to $0$. 
Else, there is $Q_1\in\N$ such that 
$$
c_\infty\in\partial(\Lpq^{H_0})
$$
for some $p/q$ with $q\leq Q_1$ and moreover
the carrots $\Dec_k$ are eventually contained inside $W_{p/q}^{H_0}$, because there are 
finitely many wakes $W_{p/q}^{H_0}$ with $q\leq Q_1$ and they are strongly separated. 

Secondly we apply Yoccoz parameter puzzles theorem, \thmref{yoccozpuzzle}. 
For any $p/q$ the corresponding rotation orbit $0<\theta_0<\ldots<\theta_{q-1}$ 
is disjoint from the set of dyadic arguments. 
Thus for any $p/q$ the graph defining the associated $p/q$ puzzle 
for $\Lpq^{H_0}$ is disjoint from $\Psi(\whDec)$. Therefore, 
there exists an increasing  sequence $n_k$ such that for $k\ge k_0$ the carrots $\Dec_k\subset \overline{\mathcal P}_{n_k}\ni c_\infty$  
(the Parameter Puzzle Piece).

Hence by Yoccoz \thmref{yoccozpuzzle} either the diameter tends to $0$ or the limiting parameter $c_\infty$ is 
renormalizable, that is $c_\infty\in \Mp$ for some period $k$ first renormalization copy $\Mp$ of $\Mbrot$ in $\Lpq^{H_0}$, 
where $q\leq Q_1$ and $q\leq k$.

\subsubsection*{The third reduction : reduction to the toy example.}
As above let $\theta_\pm$ denote the arguments of the external rays of $\Mbrot$ 
co-landing at the root $c_0'$ of $\Mp$. 
\emph{Let $\Lambda$\label{Lambdadef} denote the parameter disk whose 
closure contains $\Mp$ and which is bounded by the segments of the rays $R_{\theta_\pm}^\Mbrot$ 
with potential up to and including $2$ union a connecting subarc of the level $2$ equipotential 
$E_\Mbrot(2):=\Psi(C(0,\e^2))$,} where $C(0,\e^2) = \{z \mid \;|z| = \e^2\}$.

We need the following result on $\Mbrot$ presumably due to Douady, Hubbard and Lavaurs. 
\REFTHM{Lavaurs}
Let $0<\eta_-<\eta_+<1$ be rationals for which the parameter rays $R_{\eta_\pm}^\Mbrot$ 
coland at some point $c_0\in\Mbrot$ and let $W_{\eta_-,\eta_+}^\Mbrot$ 
denote the parameter sector bounded by $R_{\eta_-}^\Mbrot\cup\{c_0\}\cup R_{\eta_+}^\Mbrot$ 
and not containing $H_0$. 
Then the forward orbits of $\eta_\pm$ do not enter the interval $]\eta_-,\eta_+[$ 
and for any $c\in W_{\eta_-,\eta_+}^\Mbrot$ the pair of dynamical rays $R_{\eta_\pm}^c$ move 
homorphically with $c$, co-land at some 
repelling (pre)periodic point $z(c)$ with $Q_c^{k'+l}(z(c)) = Q_c^l(z(c))$, where $l\geq 0$ 
is the common preperiod of $\eta_\pm$ and $k'>0$ divides the common period 
$k>0$ of $\sigma^l(\eta_\pm)$ and the set $R_{\eta_-}^c\cup\{z(c)\}\cup R_{\eta_+}^c$ 
bounds a sector $W^c$ containing $c$, but not $0$.
\ENDTHM
\newcounter{Lavaurscounter}
\setcounter{Lavaurscounter}{\value{theorem}}
This theorem is at least folklore. But because we do not have a precise reference, we have for completeness provided a proof in the Appendix, on page~\pageref{proofofLavaurs}.

Thus for any $c\in\Lambda$ the dynamical rays $R_{\theta_\pm}^c$ co-land at a repelling $k$-periodic 
point $\beta'_c$ and the rays $R_{\whtheta_\pm}^c$ co-land at the $Q_c^k$-preimage $\cobeta'_c$ 
of $\beta'_c$ all of which moves holomorphically with $c\in\Lambda$. 
Moreover the set
$$
\bigcup_{j=0}^{k-1} Q_c^{-j}(R_{\theta_-}^c\cup\{\beta'_c\}\cup R_{\theta_+}^c) 
$$
does not enter the sector $W_{\theta_-,\theta_+}^c$ bounded by the closure of the colanding pair of rays 
$R_{\theta_-}^c$ and $R_{\theta_+}^c$ and containing $c$. Hence
$$
E_c(1)\bigcup\bigcup_{j=0}^k Q_c^{-j}(R_{\theta_-}^c\cup\{\beta'_c\}\cup R_{\theta_+}^c) 
$$
moves holomorphically with $c\in\Lambda$. 
Similarly to the relatively dyadic wakes of $\Mp$ we define the relatively dyadic wakes 
$W^c_0$ as the open set not containing $c$ and bounded by the closure of the rays $R_{\theta_\pm}^c$ and 
$W^c_{1/2}=W_{\whtheta_+,\whtheta_-}^c$ as the open set bounded by the closure of the rays $R_{\whtheta_\pm}^c$ and disjoint from $W^c_{0}$.
For any $c\in\Lambda$ there is a renormalization, a quadratic-like restriction of $Q_c^k$ 
for which the filled-in Julia set $K_c'\subset K_c$ consists of the points in the filled-in Julia set of $K_c$, 
whose orbits never enters the relatively dyadic wakes $W^c_{0}$ and $W^c_{1/2}$ 
(see also \eqref{renormalization} below). 

The key point in the proof of \propref{shrinkingcarrots} 
is that all of the dyadic carrots $\Dec_{p/2^n}$ are disjoint from $\Mp$, 
because their root points $\Psi(p/2^n)$ are disjoint from $\Mp$.
And if such a carrot intersects $\Lambda$, then it is entirely contained in $\Lambda$ 
and its dynamical counter part in the dynamical planes of $Q_c$  is either 
contained in $W^c_{1/2}$ or has a univalent forward image, which is. 
In order to prove the theorem we shall wrap the dynamical counter part of 
each dyadic carrot inside the relatively dyadic wake $W^c_{1/2}$ in an annulus in $W^c_{1/2}$ 
moving holomorphically with $c\in\Lambda$ and of modulus bounded uniformly from below.

To do this we shall follow slightly different paths according to wether $\Mp$ is 
a primitive copy or the satelite copy $\Mpq$ with root on the cardioid. 
We start with the primitive case and afterwards indicate the changes 
which make the proof in the satelite case.

\subsubsection*{The Primitive Case}
Suppose $\Mp\subset\Lpq$ is a primitive copy of $\Mbrot$. 
Let $\delta_0^c$ denote the subarc of 
$$
R_{\theta_-}\cup\{\beta'_c\}\cup R_{\theta_+} = \overline{R_{\theta_-} \cup R_{\theta_+}}
$$
consisting of points with potential up to and including $1$. 
Similarly let $\delta_{1/2}^c$ denote the subarc of $\overline{R_{\whtheta_-}\cup R_{\whtheta_+}}$
up to and including potential $1$. 
Define similarly the parameter arcs $\delta_0^\Mp$ and $\delta_{1/2}^\Mp$.
Let $c\in\Mp$ be arbitrary and let $P=P_n^c$ denote the $p/q$ puzzle piece of level $n$ 
containing $c$ given by \thmref{yoccozpuzzle}. 
Let $\eta^-<\eta^+$ denote the (rational) arguments of the co-landing pair of 
external rays for $Q_c$, which are on the boundary of $P$ and which separates $c$ from $0$. 
Then the parameter rays $R_{\eta_\pm}^\Mbrot$ co-land at some point in $\Mbrot$. 
Denote by $\Lambda^P$\label{LambdaPdef} the parameter disk which contains $\Lambda$, and 
which is bounded by a subarc of $\overline{R_{\eta_-}^\Mbrot\cup R_{\eta_+}^\Mbrot}$ union 
a subarc of $E^\Mbrot(2)$. 
Then by \thmref{Lavaurs} the dynamical rays $R_{\eta_\pm}^c$ co-land for every $c\in\Lambda_P$ and 
the arc $\overline{R_{\eta_-}^c\cup R_{\eta_+}^c}$ 
moves holomorphically with $c\in\Lambda^P$.
Denote by $\gamma_0^c$ the subarc of  $\overline{R_{\eta_-}^c\cup R_{\eta_+}^c}$ 
of potential up to and including $1$ and for $c\in\Lambda$ let $U_0^c$ denote the disk 
not containing $0$ and which is bounded by $\gamma_0^c$ 
and a subarc of the equipotential $E^c(1)$. 
Then $\partial U_0^c$ moves holomorphically over $\Lambda^P$.
Write $\Lambda_0^P :=\Lambda^P$ and $\Lambda_1^P=\Lambda^P\cap F^\Mbrot(1)$. 
Moreover for $c\in\Lambda_1^P$ let $U_1^c$ 
denote the connected component of $Q_c^{-k}(U_0^c)$ containing the ends 
$(R_{\theta_-}^c\cup R_{\theta_+}^c)\cap F^c(2^{-k})$. 
Then for $c\in\Lambda_1^P$ the restriction 
\begin{equation}\label{renormalization}
{\mapfromto {f_c := Q_c^k} {U_1^c} {U_0^c}}
\end{equation}
is quadratic like, $\partial U_1^c$ moves holomorphically with $c$ and the filled Julia set $K'_c$ is 
connected, if and only if $c\in\Mp$. 
Let $\omega_c\in U_1^c$ denote the unique critical point of $f_c$, 
so that $f_c(\omega_c) = Q_c^k(\omega_c) = c$. 
Notice that $\Lambda_1^P$ is precisely the set of parameters for which $c\in U_0^c$, 
in fact $c\in \partial U_0^c$ if and only if $c\in \partial\Lambda_1^P$.
For later use we define $U_n^c = f_c^{-n}(U_0^c)$, which may or may not be connected for 
$n>1$. 

\begin{figure}[htbp]
\centerline{\includegraphics[height=65mm]{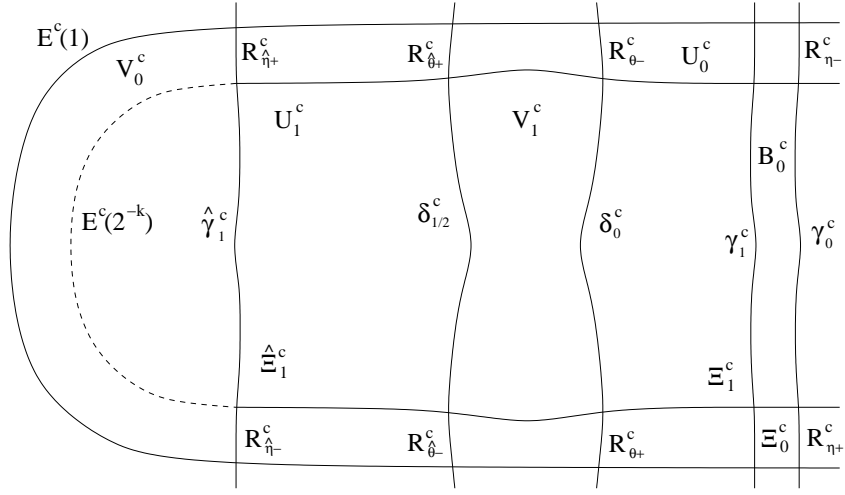}}
\caption{The disks $U_0^c$, $U_1^c$, and $\whXi_1^c$. 
The set $U_0^c$ is the disk insided $E^c(1)$ and to the left of $\gamma_0^c$. 
The set $U_1^c$ is the disk inside $E^c(2^{-k})$ and bounded to the right and left 
by $\gamma_1^c$ and $\whgamma_1^c$ respectively. 
The arc $\delta_0^c$ separates the two subdisks $V_0^c$ and $\Xi_0^c$ 
of $U_0^c$, $V_0^c$ to the left and $\Xi_0^c$ to the right of $\delta_0^c$.
The subsets $\whXi_1^c$, $V_1^c$ and $\Xi_1^c$ of $U_ 1^c$ are to 
the left of $\delta_{1/2}^c$, between $\delta_{1/2}^c$ and $\delta_0$ and to the right of $\delta_0$
repsectively.
\label{basicdisposition}}
\end{figure}
\nobreak
Write $\Lambda_0=\Lambda\subset\Lambda_0^P$ and $\Lambda_1=\Lambda\cap\Lambda_1^P$.
For $c\in\Lambda_0$ let $V_0^c$ be the connected component of $U_0^c\sm\delta_0^c$ 
containing $\omega_c$ let $\Xi_0^c$ denote the other connected component. 
Define recursively $V_n^c = f_c^{-n}(V_0^c)$ and $\Xi_n^c=f_c^{-n}(\Xi_0^c)\cap\Xi_0^c$ 
(see also \figref{basicdisposition}).
Then the restriction {\mapfromto {f_c} {V_{n+1}^c}{V_n^c}} is a $2:1$ branched covering, 
whereas {\mapfromto {f_c} {\Xi_{n+1}^c}{\Xi_n^c}} is an isomorphism.

Let $\gamma_1^c$ denote the extension to potential level $1$ of 
$Q_c^{-k}(\gamma_0^c)\cap\partial\Xi_0^c$ and let $B_0^c\subset \Xi_0^c$ denote the quadrilateral 
bounded by $\gamma_0^c$, $\gamma_1^c$ and subarcs of $E^c(1)$. 
Define recursively the univalently iterated preinages $B_{n+1}^c = f_c^{-1}(B_n^c)\cap\Xi_n^c$. 

For each $n\geq 1$ let $\whXi_n^c$ denote the ``other'' connected component of $f_c^{-1}(\Xi_{n-1}^c)$, 
having a boundary arc in $\delta_{1/2}^c$. 
For each $n\geq 1$ let  $\whB_n^c\subset\whXi_n^c$ 
denote the ``twin'' of $B_{n}^c$, i.e. the connected component of $f_c^{-1}(B_{n-1}^c)\cap\whXi_n^c$. 
Let $\whgamma_1=Q_c^{-k}(\gamma_0^c)\cap\partial\whXi_1^c$ 
extended to equipotential level $3/4$ 
and let $\Omega^c$ denote the open disk bounded by $\whgamma_1^c$ 
and the subarc of $E^c(3/4)\cap U_0^c$ connecting the endpoints of $\whgamma_1^c$. 
Let $\Dc{}\subset V_0^c$, denote the disc bounded by $\delta_{1/2}^c$ 
union the subarc of $E^c(1)$ connecting the endpoints of $\delta_{1/2}^c$. 
Then by construction each of the sets 
$\Omega^c$ and $\whB_n^c$, $n\geq 1$ are relatively compact in $\Dc{}$. 
In fact 
\REFLEM{uniformbound}
For every $c\in\Lambda$ there exists $m= m(c) >0$ such that 
$$
\mod(\Dc{}\sm\Omega^c)\geq m,\quad\textrm{and}\quad\forall\;n\geq 1\quad 
\mod(\Dc{}\sm\overline{\whB_n^c})\geq m
$$
so that
$$\forall\;n\geq 1\quad 
\mod(\Xi_0^c\sm\overline{B_n^c})\geq m.
$$
\ENDLEM
\begin{figure}[htbp]
\centerline{\includegraphics[height=65mm]{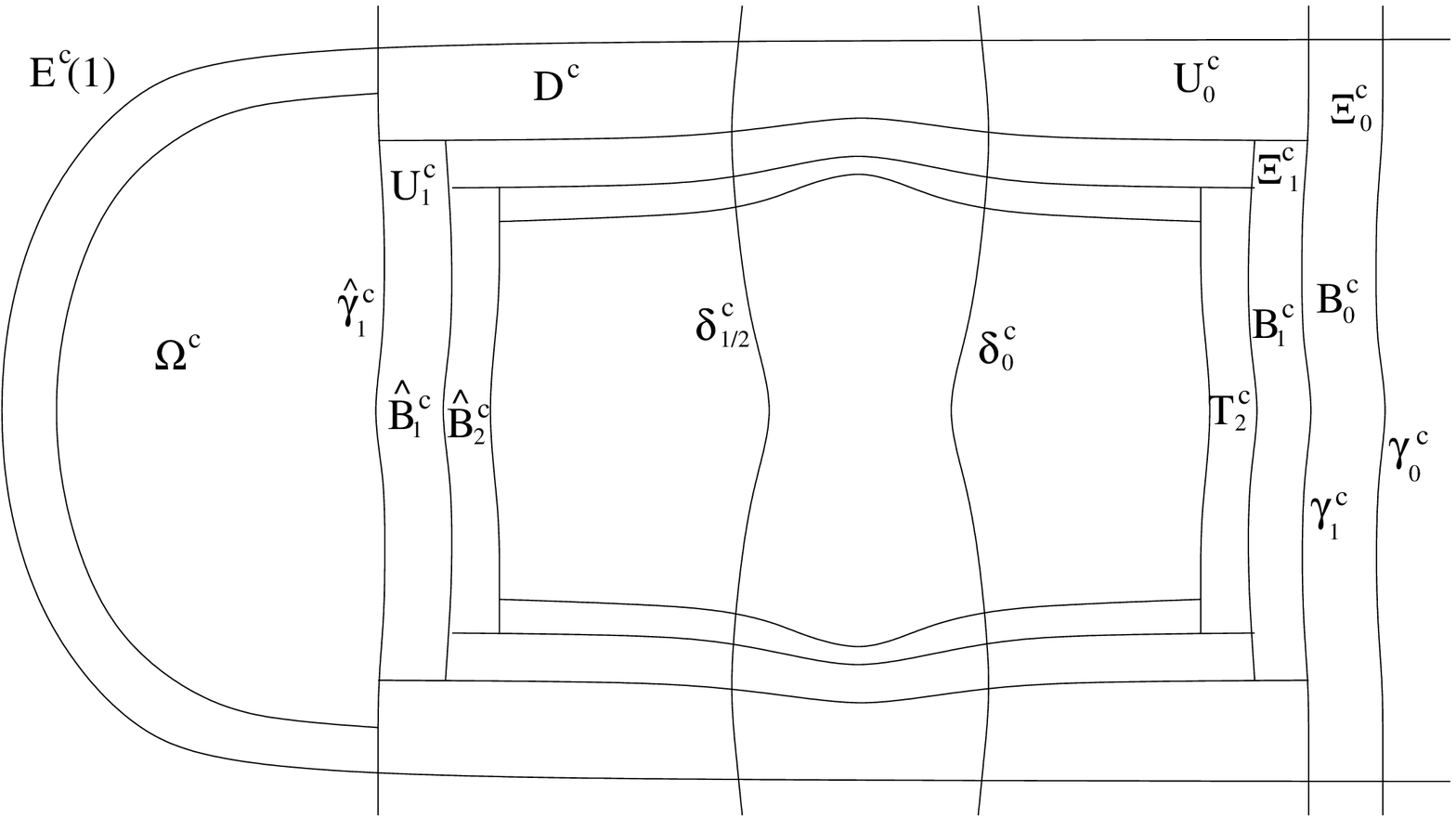}}
\caption{The decomposition of the disk $\Dc{}$, which is bounded by the equipotential $E^c(1)$ and 
the arc $\delta_{1/2}^c$. The disk $\Xi_0^c$ is the subset of $U_0^c$ 
to the right of $\delta_0^c$. The disk $\whXi_1^c$ (not labelled) is the subset of $U_1^c$ 
to the left of $\delta_{1/2}^c$. The disc $\Omega^c$ is to the left of $\whgamma_1^c$ 
and inside the equipotential $E^c(3/4)$ (not labelled).\label{boxdisposition}}
\end{figure}

\PROOF
The restriction {\mapfromto {f_c} {\whXi_1^c} {\Xi_0^c}} is biholomorphic 
so that for all $n\geq 1$
\begin{align*}
\mod(\Dc{}\sm\overline{\whB_{n+1}^c})&\geq\mod(\whXi_1^c\sm\overline{\whB_{n+1}^c}) = \mod(\Xi_0^c\sm\overline{B_n^c})\\
&\geq\mod(\Xi_{n-1}^c\sm\overline{B_n^c}) 
= \mod(\Xi_0^c\sm\overline{B_1^c})
 \end{align*}
Thus we may define 
$$
m(c) = \min\{\mod(\Dc{}\sm\Omega^c), \mod(\Dc{}\sm\overline{\whB_1^c}), \mod(\Xi_0^c\sm\overline{B_1^c})\}.
$$
\ENDPROOF

Moreover again by construction the graph
$$
G^c =\partial \Dc{}\cup\partial\Omega^c\cup\bigcup_{n=1}^\infty (\partial\whXi_n^c\cup\partial\whB_n^c)\cup 
\bigcup_{n=0}^\infty (\partial\Xi_n^c\cup\partial B_n^c) \cup\partial U_0^c
$$
moves holomorphically with $c\in\Lambda_0$ via the {\Bottcher}-coordinates at $\infty$, 
and does not intersect any of the dynamical dyadic carrots $\Dec_{p/2^m}^c$. 
The later because all such carrots are at dyadic angles and only $\Dec_0^c$ extends 
further than potential $1/2$. 
The holomorphic motion part because,  $Q_c^n(G^c)$ does not meet the critical point $0$ 
for any $n\geq 0$ or $c\in\Lambda$:
When the parameter $c\in\Lambda_0=\Lambda$ the critical value $c$ does not belong to 
$\overline{\Xi_0^c}$. 
Hence the {\Bottcher}-coordinate is defined and depends holomorphically on $c$, 
on the dense subset $(\partial\Xi_0^c\cup\partial B_0^c)\sm K_c$ 
of $(\partial\Xi_0^c\cup\partial B_0^c)$, so that the later moves holomorphically with $c$. 
Secondly $f_c$ depends holomorphically on $c$ and its critical value $c$ again 
still does not belong to $\overline{\Xi_0^c}$. 
Hence the iterated univalent preimages 
of $(\partial\Xi_0^c\cup\partial B_0^c)$ inside $\Xi_0^c$ and $\overline{\whXi_1^c}$ 
depend also holomorphically on $c$. 
This takes care of 
$$\bigcup_{n=1}^\infty (\partial\whXi_n^c\cup\partial\whB_n^c)\cup 
\bigcup_{n=0}^\infty (\partial\Xi_n^c\cup\partial B_n^c).
$$
Finally $\partial\Dc{}$ moves holomorphically with $c$, because $\delta_{1/2}^c$ does and 
$\partial\Omega^c$ moves holomorphically with $c$, because $\whgamma_1^c$ does.

No dynamical plane dyadic carrot $\Dec_{p/2^n}^c$ intersects $G^ c$.  
Hence any such dyadic carrot in the relative dyadic wake $\Wc_{1/2}$ of 
the filled Julia set $K'_c$ for {\mapfromto {f_c} {U_1^c} {U_0^c}} 
is contained in one of the sets $\Omega^c$ or $\whB_n^c$ for some $n\geq 1$. 
Such a carrot is thus wrapped by an annulus of modulus at least $m(c)>0$, 
contained in $\Dc{}$ and thus disjoint from $K'_c$. 

We shall use an argument to transfer bounds on moduli of dynamical annuli 
to bounds on moduli of corresponding parameter annuli. 
This argument was pioneered by Shishikura (see also \cite{Roesch}).

Define $\Lambda_1=\Lambda\cap F(1)$ and fix as basepoint $c_b\in\Lambda_1\subset\Lambda_0$ 
the center of the central hyperbolic component of $\Mp$. 
Note that on $F(2)\supset \Lambda_0$ the {\Bottcher}-coordinate at $\infty$ defines a 
holomorphic motion of the set $\Cbar\sm F^{c_b}(1)$ extending the {\Bottcher}-motion of $\Gcb$.
By Slodkowski's extension theorem there exists a global holomorphic motion 
{\mapfromto H {\Lambda_0\times\Cbar}\Cbar} over $\Lambda_0$ 
with base point $c_b$ and extending the {\Bottcher}-motion of the graph $\Gcb$ union $\Cbar\sm F^{c_b}(1)$, 
in particular we obtain holomorphic motions of $\overline{U_0^{c_b}}\supset\Dcbbar{}$.  
As usual for $c\in\Lambda$ write $H_c(\cdot)=H(c,\cdot)$, then 
each map {\mapfromto {H_c} {\Cbar} {\Cbar}} is a quasi-conformal homeomorphism 
with a dilatation bounded uniformly above by $\log\d_\Lambda(c,c_b)$, 
where $\d_\Lambda(\cdot,\cdot)$ denotes hyperbolic distance in $\Lambda=\Lambda_0$. 

Define similarly to $\Dc{}$ the parameter disk $\DMp\subset \Lambda_1\sm\Mp$ 
as the disc bounded by $\delta_{1/2}^\Mp$ union the subarc of $E^\Mbrot(1)$ 
connecting the endpoints of $\delta_{1/2}^\Mp$. 
Then $\DMp{}$ is relatively compact in $\Lambda_0$.
Let {\mapfromto {\chi} {\Lambda_0} \Cbar} be the map $\chi(c) = H_c^{-1}(c)$. 
Then $\chi$ is locally quasi-regular with dilatation $K(c)$ bounded by $\log\d_\Lambda(c,c_b)$. 
By construction the restriction {\mapfromto {\chi} {\partial\DMp} {\partial\Dcb{}}} 
is a homeomorphism. 
Hence the restriction {\mapfromto {\chi} {\DMpbar{}} {\Dcbbar{}}} is the restriction of 
a quasi-conformal homeomorphism with dilatation bounded by 
$K=K_{1/2} = \max\{\log\d_\Lambda(c,c_b)| c\in\DMpbar\}$.

Define $\Omega^\Mp = \chi^{-1}(\Omega^{c_b})\subset\DMp{}$ and 
$\whBMp_n=\chi^{-1}(\whBMp_n)\subset\DMp{}$, $n\geq 1$. 
Then any dyadic carrot $\Dec_{p/2^n}\subset W'_{1/2}$ is contained in one of 
the disks $\Omega^\Mp$ or $\whBMp_n$ and is thus wrapped in an annulus with a modulus 
bounded from below by $K_{1/2}\cdot m(c_b)$ according to \lemref{uniformbound}.

Rename $\Dc{} =: \Dc{1/2}$, $\DMp{} =: \DMp{1/2}$ and 
$\chi=:\chi_{1/2}$. 
We have proved that any dyadic carrot $\Dec_{p/ 2^n}$ in the relative $1/2$ wake 
$W'_ {1/2}$ of $\Mp$ is wrapped in an annulus of modulus uniformly bounded from below and 
contained in the disk $\DMp{1/2}$, which is disjoint from $\Mp$. 
Moreover the annuli are q-c images of corresponding annuli in the dynamical plane of $Q_{c_b}$.
We shall prove by induction the similar statements for any other relative dyadic wake $W'_{r/2^s}$ 
of $\Mp$. The only difference is that the bounds on the dilatation of the q-c homeomorphisms 
$\chi_{r/2^s}$ and hence on the moduli of annuli in the $W'_{r/2^s}$ wake depends on $r/2^s$. 
As a remedy for this we shall apply the Levin-Yoccoz parameter inequality once more. 
Here follow the details.

Recall that $\Vcz$ is the connected component of $U_0^c\sm\delta_0^c$ 
containing the critical point $\omega_c$ and $\Vcn=f_c^{-n}(\Vcz)$. 
We shall need also the extension $\wtVco=\Vcz\sm\Dcbar{}$ of $\Vco$ and its iterated 
preimages $\wtVcn=f_c^{-n}(\wtVco)$.
Define parameter disks  $\Lambda_s$, $s>1$ by
$$
\Lambda_s = \{c| c \in \wtVc_{s-1}\}.
$$
Evidently $\Lambda_s\supset\Lambda_ {s+1}$.
Note that the condition $c\in\wtVc_s$ is equivalent to $f_c^s(\omega_c)\in\wtVc_1$.
Rename $\Gc =: \Gc_1$ and define recursively , 
$\Gc_{s+1}= f_c^{-1}(\Gc_s)\cup\Gc_1$ for $s\geq 1$. 
Then as noted above $G_1^c$ moves holomorphically in $\Lambda_0\supset\Lambda_1$ 
and we shall prove as part of the induction on $s\geq 2$, 
that for $c\in\Lambda_s$ the critical value $c\notin\Gc_{s-1}$, 
so that $\Gc_s$ moves holomorphically over $\Lambda_s$. 

For $s=2$ notice that, by the above $c\in\Dcbar{1/2}$ if and only if $c\in\DMpbar{1/2}$. 
Thus $c\notin\Gc_1$ for any $c\in\Lambda_2$, so that 
$G_2^c$ moves holomorphically with $c\in\Lambda_2$. 
For $c\in\Lambda_2$ let $\Dc{rj/2^2}$ for $r=1,3$ denote the connected components of $f_c^{-1}(\Dc{1/2})$ 
containing the $r/2^2$ dyadic carrots and define $\DMp{r/2^2}\subset\Lambda_2$ as the 
parameter disks bounded by the corresponding parameter ray segments and equipotential level. 

Rename the previous holomorphic motion $H$ to $H_1$ and let $H_2$ denote the restriction of $H_1$ 
to $\Lambda_2\times(\Cbar\sm\wtVcb_1)$. Extend $H_2$ to a motion including $\wtVcb_1\sm\Vcb_1$ 
using the {\Bottcher}-motion and extend $H_2$ further to 
$\Lambda_2\times f_{c_b}^{-1}(\Dcbbar{1/2})$ 
by $f_c^{-1}\circ H_1(c,f_{c_b}(z))$,
where the inverse branches are taken so as to map $\Dcb{r/2^2}$ quasi-conformically onto $\Dc{r/2^2}$. 
Finally use Slodkowski's extension theorem to extend this holomorphic motion to a holomorphic motion 
of $\Cbar$ over the disk $\Lambda_2$ (i.e. extend the motion by a motion of $\Vcb_2$). 
By the same argument as above the map {\mapfromto {\chi_2} {\Lambda_2} {\wtVcb_1}} 
given by $\chi_2(c) = {(H_2)_c^{-1}(c)}$ is a locally quasi regular map. 
Again by construction {\mapfromto {\chi_2} {\partial\DMp{r/2^2}} {\partial\Dcb{r/2^2}}} 
are homeomorphisms so that the restrictions {\mapfromto {\chi_2} {\DMp{r/2^2}} {\partial\Dcb{r/2^2}}} 
are quasi-conformal. 
However on the sets $\DMp{r/2^2}$ the holomorphic motion $H_2$ is a conjugacy between 
the holomorphic maps $f_{c_b}$ and $f_c$. Hence the dilatation of ${(H_2)}_c$ at $z$ equals that of 
${(H_1)}_c$ at $f_{c_b}(z)$. 
Hence again the dilatation of $\chi_2$ on $\DMp{r/2^2}$ is again bounded by 
bound given by $K_{r/2^s} = \max\{\log\d_\Lambda(c,c_b)| c\in\DMpbar{r/2^2}\}$. 
Arguing as in the initial case corresponding to $s=1$ completes the case $s=2$. 
The inductive step is completely similar and is left to the reader.

Let $H'$ denote the central hyperbolic component of $\Mp$. 
Then for $k$ sufficiently large $c_k$ belongs to the $p'_k/q'_k$ limb $L_{p'_k/q'_k}^{H'}$ 
of $H'$. 
Applying the Yoccoz-Levin parameter inequality \thmref{yoccozinequality} to $H'$ 
we find that the diameter of $L_{p'/q'}^{H'}$ union its attached dyadic carrots is 
bounded uniformly by $C/q'$ for some constant $C=C_{H'}$. 
Thus it suffices to consider the case $q'_k\leq Q'$ for some integer $Q'$. 
As the containing wakes $W_{p'/q'}^{H'}$, $q'\leq Q'$ are strongly separated we can 
further assume that $p'_k/q'_k = p'/q'$ for $k$ large. 
The set  $W_{p'/q'}^{H'}\cap\Lambda_1$ is relatively compact in $\Lambda$
so that 
$$
\sup\{\log\d_\Lambda(c,c_b)| c\in W_{p'/q'}^{H'}\cap\Lambda_1\} = K= K_{p'/q'}^{H'} < \infty.
$$
Hence the dyadic carrots $\Dec_k$ either has a diamter which a priori tends to zero or such 
carrots are separated from $\Mp$ by an annulus in $\Lambda\sm\Mp$ of modulus at least $m(c_b)/K$. 
And in the latter case their diameters are forced to converge to zero a posteriori. 
Becuase the roots $c_k\in\Dec_k$ converge to $c_\infty\in\Mp$, 

This completes the proof that if $\Lambda\ni c_k \to c_\infty\in\Mp$, then 
the diameter of $\Dec_k$ converge to zero. 
For the case $c_k\in\Lambda_P\sm\Lambda_0$ we necessarily have 
$c_\infty = c_r$, where $c_r$ denotes the root of $\Mp$. 
To prove that the diameter of $\Dec_k$ converge to zero also in this case let
$$
\Lambda_n^P = \{c\in\Lambda_P| c\in U_n^c\}.
$$
For any $c\in\Lambda=\Lambda_0$ the sets $\partial\Xi_n^c$, $n\geq 0$ move holomorphically with $c$. 
Define $A_n^c = \Xi_n^c\sm\Xi_{n+1}^c$, then the $A_n^c$ are quadrilaterals with a-sides 
the boundary arcs  $\partial A_n^c\cap R_{\theta_-}^c$ and $\partial A_n^c\cap R_{\theta_+}^c$. 
Moreover  $\partial A_n^c$ even move holomorphically with $c\in(\Lambda_0\cup\Lambda_n^P)$. 
Let 
$$
A_n^\Mp = \{c\in\Lambda_n^P| c \in A_n^c\}
$$
denote the corresponding parameter quadrilaterals. 
Then the root $c_r$ of $\Mp$ belong to $\Lambda_n^P$, for all $n$. 
Choose by Slodkowski's extension theorem a holomorphic motion 
$$
H^0:\Lambda_0^P\times A_0^{c_r} \longrightarrow \C
$$
over $\Lambda_0^P$ with base point $c_r$ of the quadrilateral 
$A_0^{c_r}$ extending the B{\"o}ttcher motion of its boundary. 

For $c\in\Lambda_n^P$ the restriction {\mapfromto {f_c^n} {A_n^c} {A_n^0}} is biholomorphic. 
Hence we may lift the motion $H^0$ to a holomorphic motion 
$$
H^n:\Lambda_n^P\times A_n^{c_r} \longrightarrow \C. 
$$
As with the annuli above define quasi-conformal homeomorphisms
$$
\rho_n: A_n^\Mp \longrightarrow A_n^{c_r}
$$
by $\rho_n(c) = {\left(H_c^n\right)}^{-1}(c)$. 
Then as above these have q.-c.~distortion bounded by the distortion of the q.-c.~homeomorphisms 
$H_c^0(\cdot)$, $c\in A_n^\Mp$. That is bounded by
$$
K = \sup\{\log\d_{\Lambda_0^P}(c,c_r)| c\in A_n^\Mp\}
$$
which is uniformly bounded, because $A_n^\Mp\subset\Lambda_1^P\subset\subset\Lambda_0^P$. 
Thus all the quadrilaterals $A_n^\Mp$ have modulus bounded uniformly from below
by $\mod(A_0^{c_r})/K$. Moreover the a-sides of these quadrilaterals are all contained 
in the two rays $R_{\theta_-}^\Mp$ and $R_{\theta_+}^\Mp$ co-landing at $c_r$. 
By the Gr{\"o}tzsch-inequality for annuli the euclidean diameter of $A_n^\Mp$ tend 
to zero and the closures converge to $c_r$. By construction no dyadic carrot intersects 
the boundary of any of the $A_n^\Mp$. 
Thus also in this case the diameter of $\Dec_k$ converge to zero as $k\to\infty$.
This completes the proof in the case $c_\infty$ belongs to a primitive first renormalization copy.

\subsubsection*{The satelite case}
In the complementary satelite case $\Mp=\Mpq$ with central hyperbolic component $\Hpq$ 
attached at internal argument $\exp(i2\pi p/q)$ from the central hyperbolic component $H_0$ of $\Mbrot$. 
Let as above $\theta_-<\theta_+$ be the arguments of the parameter rays co-landing at the root and bounding 
the wake $W_{p/q}^{H_0}$. 
Recall that $c_\infty\in\Mpq$ is the limiting parameter of the roots of 
dyadic carrots 
and that these 
dyadic carrots are eventually 
contained in $W_{p/q}^{H_0}$.

We apply the Yoccoz-Levin parameter inequality \corref{yoccozinequality} similarly as we 
have done twice above. 
This reduces the problem to the case where $c_\infty$ belongs to 
a relative $p'/q'$-limb $L_{p'/q'}^{\Hpq}$ of $\Hpq$ for some $p'/q'\not=0/1$ 
and the dyadic carrots $\Dec_k$ are subsets of the corresponding wake 
$W_{p'/q'}^{\Hpq}$ for large $k$. 
Denote by $\tau_-<\tau_+$ the arguments of the parameter rays bounding $W_{p'/q'}^{\Hpq}$ and 
define $\Lambda=\Lambda_0 = W_{p'/q'}^{\Hpq}\cap F^\Mbrot(2)$. 
Let $\IMp = I^\Mpq$ denote the Cantor set of arguments of parameter rays accumulating $\Mp$ 
as given by the Douady tuning algorithm. 
Then $\tau_\pm\in\IMp$ and each has a unique preimage $\whtau_\pm=\sigma^{-q}(\tau_\pm)\cap\IMp$ 
different from itself. 
For $c\in\Lambda$ let $U_0^c$ denote the disk containg the fixed point $\alpha^c$ of $Q_c$ and bounded 
by the segments $\overline{(R_{\sigma^i(\whtau_-)}^c\cup R_{\sigma^i(\whtau_+)}^c)\cap F^c(1)}$ for 
$0<i<q$ union the connecting subarcs of $E^c(1)$. 
Denote by $\iota^c$ the open subarc of $\partial U_0^c\cap\Ec(1)$ intersecting the rays $R_{\theta_\pm}^c$ and 
let $\gamma_0^c=\partial U_0^c\sm\iota^c$. 
As in the primitive case write $\delta_0^c=\overline{(R_{\theta_-}^c\cup R_{\theta_+}^c)\cap F^c(1)}$ 
for $c\in\Lambda$. 

Then the whole setup is similar to the primitive case. 
 We can thus define $\Omega^c$, $\Xi^c_n, \Bc_n, \whXi^c_{n+1}, \whBc_{n+1}$ for $n\geq 0$ and 
 $\Gc$, all of which moves holomorphically with $c\in\Lambda$.
There are however two differences:
The first is that the center of $\Hpq$ does not belong to $\Lambda$. 
The arguments we used in the primitive case are in-sensitive to a change of base point 
$c_b$ to another point in the interior of $\Mbrot$. We shall thus take as base point $c_b\in\Lambda$ 
the center of the central hyperbolic component $H_{p'/q'}^{\Hpq}\subset W_{p'/q'}^\Hpq$. 
The second difference is that the Yoccoz-Levin parameter inequality is 
applied to the sublimbs of the hyperbolic component $H_{p'/q'}^{\Hpq}$.
We leave the details to the reader.

This completes the satelite case and thus completes the proof of \propref{shrinkingofconvgcarrots}.


\subsubsection*{Proving \thmref{shrinkingcarrotsp}}
As external rays do not cross the proof of \thmref{shrinkingcarrotsp} is completely analogous 
to the proof above of \thmref{shrinkingcarrots}. 
Let $\Mp$ be any copy of $\Mbrot$ inside $\Mbrot$ or $\Mone$. 
In the arguments above replace the carrot field $\Dec$ of $\Mbrot$ 
by the dyadic decorations $\Decp$ of $\Mp$. 
Use Yoccoz parameter inequality and the iterated Yoccoz parameter puzzle theorem 
relative to $\Mp$ to prove that: $\diam(\Decp_k)\to 0$ for any sequence of decorations $(\Decp_k) _k$ 
with roots $c_k$ converging to a relatively non renormalizable parameter $c_\infty\in\Mp$. 
Secondly consider the case $c_\infty\in M''$, where $M''\subset \Mp$ 
is a relative to $\Mp$ first renormalizable copy of $\Mbrot$ belonging to some $p/q$ limb 
of the central hyperbolic component $H'$ of $\Mp$. 
Use again the Yoccoz puzzle relative to $\Mp$ to define the parameter disk $\Lambda$ 
containing $M''$ similarly as we defined $\Lambda$ for $\Mp$ above. 
And define also $\Lambda^P$ analogously, i.e.~with the aid of the $p/q$ puzzle piece $P$ relative to $\Mp$ 
given by \thmref{yoccozpuzzle} for $\Mp$. From here the proof proceeds analogously.

\section*{Appendix}
In this appendix we supply for completeness proofs of the two theorems we refer to, 
but for which we have not been able to find either adequate or complete proofs in the litterature.

\setcounter{theorem}{\value{YoccozLevincounter}}
\addtocounter{theorem}{-1}
\THM[The Levin-Yoccoz Dynamical Inequality]
Let $H$ be any hyperbolic component of $\Mbrot$ of period $k$. 
Let $p/q$ be any non zero reduced rational and let $W_{p/q}^H$ denote 
the relative $p/q$ wake of $H$, bounded by parameter rays with arguments 
$0<\eta_-<\eta_+<1$. 
For any $c\in W_{p/q}^H$ let $\lambda$ denote the multiplier of the repelling 
$k$-periodic common landing point $\alpha'$ of the $kq$ periodic rays $R_{\eta_\pm}^c$. 
Then $\alpha'$ has combinatorial rotation number $p/q$ and 
 $\lambda$ has a logarithm $\Lambda$ such that:
$$
|\Lambda -p/q2\pi i| \leq \frac{2k\log2\cos\theta}{q}\frac{\pi}{\omega(c)}, 
$$
where $\theta\in\;]-\pi/2,\pi/2[$ is the argument of $\Lambda -p/q2\pi i$ 
and $\omega(c)$ is the angle of vision of the interval $i2\pi[\eta_-,\eta_+]$ 
from $\Log \phi_c(c)\in\{z=x+iy| 0< y < 2\pi\}$.
\ENDTHM
\PROOF\label{YoccozLevindynamicalproof}
Levin proved the fixed point case $k=1$ in \cite[TH. 5.1]{Levin}, 
the general case is similar. For completeness we sketch here a proof. 
Let us first recall the proof of Yoccoz inequality (or the Pommerenke-Levin-Yoccoz inequaltiy), 
full details can be found in \cite{Petersen1}. 
Let $T$ denote the quotient torus $T=D^*/Q_c^{k}$, where $D^*=\{z|0<|z-\alpha'|<r\}$ and 
$r>0$ is chosen so small that $Q_c^k$ is univalent on $D=D^*\cup\{\alpha'\}$ and 
$D\subset\subset Q_c^k(D)$. Let {\mapfromto {\Pi} {D^*} T} denote the natural projection.
The two rays $R_{\eta_\pm}^c$ belong to the same orbit 
and define combinatorial rotation number $p/q$ for $\alpha'$. 
Let $\gamma=\Pi(D^*\cap R_{\eta_-}^c) = \Pi(D^*\cap R_{\eta_+}^c)$. 
Then $\gamma$ is a Jordan curve and thus the pair $(T,\gamma)$ has a conformal modulus 
which satisfies a Gr{\"o}tzsch inequality. 

Let $w_\pm=\exp(i2\pi \eta_\pm)$. 
Then $Q_0^k(w_-) = w_+$ and $Q_c^{kq}(w_\pm) = w_\pm$. 
Similarly to $T$ above let $\whT$ denote the quotient torus $\whT=\whD^*/Q_c^{k}$, 
where $\whD^*$ is a small punctured disk centered at say $w_-$ and let 
{\mapfromto {\whPi} {\whD^*} \whT} denote the natural projection. 
Then $\whPi(\whD^*\cap\Sen)$ are two disjoint Jordan curves in $\whT$, 
with complement two disjoint, symmetric and straight annuli $A_i$ and $A_o$. 
Moreover $\whgamma=\whPi(\whD^*\cap R_{\eta_-})$ is the Jordan equator of $A_o$ and 
$$
\mod(A_i) + \mod(A_o) = 2 \mod(A_o) = \mod(\whT,\whgamma).
$$
If $c\in\Mbrot$ so that $K_c$ is connected, the B{\"o}ttcher coordinate at $\infty$ 
induces an isomorphism between $A_o$ and $\Pi(S)$ where $S$ is the connected 
component of $D^*\cap B^c(\infty)$ containing the end of $R_{\eta_-}^c$. 
Hence the Gr{\"o}tzsch inequality for $(T,\gamma)$ implies that 
\begin{equation}\label{Grotzsch}
\mod(A_o) \leq \mod(T,\gamma).
\end{equation}
Writting out the values of these two numbers explicitly yields the Yoccoz dynamical inequality: 
The torus $T$ is isomorphic to $\Cstar/\lambda z$ via the linearizer for $Q_c^k$ at $\alpha'$, 
or equivalently to $\C/(\Z\Lambda+\Z i2\pi)$ via the log-linearizer. 
Let {\mapfromto {\Pi_u} \C T} denote the universal covering corresponding to the latter isomorphism. 
Then the Jordan curve $\gamma=\Pi(D^*\cap R_{\eta_-}^c)$ lifts under $\Pi_u$ 
to an arc $\Gamma$, which is invariant under the translation $z\mapsto z+ L$, 
where $L=q\Lambda-p i2\pi$ for some appropriate logarithm $\Lambda$ of $\lambda$. 
A simple computation shows that 
$$
\mod(T,\gamma) = \frac{2\pi\cos \theta}{q|L|}
$$
where $\theta$ is the angle between the vector $L$ and the positive real axis.
A similar computation shows that 
$$
2\mod(A_o) = \mod(\whT,\whgamma) = \frac{2\pi}{kq\log 2}.
$$
Hence \eqref{Grotzsch} is equivalent to
\begin{equation}\label{yoccozdyninequality}
|\Lambda-\frac{p}{q}i2\pi| \leq \frac{2k\log 2 \cos\theta}{q},
\end{equation}
which is Yoccoz inequality.

If $c\notin\Mbrot$ let $0\leq\theta<1$ denote the argument of $c$, i.e. $c\in R_\theta^c$. 
Then the B{\"o}ttcher coordinate $\phi_c$ at infinity does not extend 
to a biholomorphic map between  $B^c(\infty)$ and $\Cbar\sm\Dbar$, but almost:
It extends to a biholomorphic map of $\Cbar\sm F^c(h)$ onto $\Cbar\sm \Dbar(\e^h)$ 
where $h=g_c(c)/2$. 
Let $\psi_c$ denote the inverse of this extension, then $\psi_c$ extends continuously to 
$C(0,\e^h)$, but this extension is not injective because 
$0 = \psi_c(\exp(h+2\pi i\theta/2))=\psi_c(\exp(h+i2\pi(\theta+1)/2))$. 
Let $C=g_c(c)+i2\pi\theta$, $N_0 = [\e^{i2\pi\theta},\phi_c(c)]$ and $N_n =Q_0^{-n}(N_0)$. 
Define $N_\theta^0 = \cup_{n\geq 0} N_n$ and $N_\theta^1 = \cup_{n\geq 0} N_1$. 
Then $Q_0(N_\theta^1) = N_\theta^0$ and 
$\psi_c$ is easily seen to extend by iterated lifting to a univalent map 
from $\Cbar_\theta:=\Cbar\sm(\Dbar\cup N_\theta^1)$ into $B^c(\infty)$. 
The map $Q_0$ lifts under $\exp(z)$ to the map $z\mapsto 2z$ on $\C$. 
That is $\exp$ is a simultanuous linearizer for all the repelling periodic points of $Q_0$.
The corresponding lifted sets $\wtN_\theta^j = \log(N_\theta^j)$, $j=0,1$  are invariant 
under translation by $i2\pi$ and $2\wtN_\theta^1=2\wtN_\theta^0$. 
Thus if $w = \exp(i2\pi \tau)\in\Sen$ is periodic and if $0\leq\tau<1$ 
does not belong to the orbit of $\theta$, 
then $\C_\theta$ contains a definite sector around the horizontal 
$\wtR_\tau= \{t+i2\pi\tau| t> 0\}$, which projects to $R_\tau^0$ under $\exp$: 
Let $\tau_l<\tau<\tau_r$ be the arguments closest to $\theta$ of points in the orbit of $w$. 
Then the sectors $\wtS_l$ between $\wtR_{\tau_l}= \{t+i2\pi\tau_l| t> 0\}$ and the oblique 
line through $i2\pi\tau_l$ in the direction $v_l=C-i2\pi\tau_l$, 
and $\wtS_r$ between $\wtR_{\tau_r}= \{t+i2\pi\tau_r| t> 0\}$ and the oblique 
line through $i2\pi\tau_r$ in the direction $v_r=C-i2\pi\tau_r$ are contained in $\C_\theta$:
If not some line segment $L$ with $\exp(L)\in N_n$ for some $n\geq 1$ intersects 
say $\wtS_l$. 
But then $2^nL$ intersects the sector $2^n\wtS_l$ with top point $2^n\tau_l$,
and is also congruent modulo $i2\pi$ to $L_0=[i2\pi\theta,C]$ with $\exp(L_0) = N_0$.
Since the $2^n\tau_l$ is an argument for a point in the orbit of $w$ 
this contradicts that $\tau_l$ is the closest such argument for points in the orbit of $w$.

Consequently the sector $\wtS$ around $\wtR_\tau$ bounded by the two lines 
through $i2\pi\tau$ and of directions $v_l$ and $v_r$ is contained in $\C_\theta$.

In the case at hand $c\in W_{p/q}^H$ implies that $\eta_<\theta<\eta_+$ 
and for $\eta=\eta_-$ we have $\eta_-=\eta_l$, $\eta_+=\eta_r$.
Let $\omega_l$ and $\omega_r$ denote the angle of inclination 
of the vectors $C-i2\pi\tau_l$ and $C-i2\pi\tau_r$ respectively. 
Then the opening angle $\omega$ of $\wtS$ equals $\omega_r-\omega_l$ 
and the sector $\wtS$ projects to a straight subannulus $A_o^\theta$ of $A_o$ with 
$$
\mod(A_o^\theta) = \frac{\omega}{\pi} \mod(A_o).
$$
Arguing as for the proof of the Yoccoz inequality we obtain 
$$
\mod(A_o^\theta) = \frac{\omega}{\pi} \mod(A_o) \leq \mod(T,\gamma).
$$
Properly rewritten as with \eqref{yoccozdyninequality} above, 
this is the Levin-Yoccoz inequality except for the interpretation of $\omega$. 
This interpretation is however an elementary exercise in planar geometry and is left to the reader. 
By continuity the inequality even holds on $\partial W_{p/q}^H$, where either 
$\omega_l$ or $\omega_r$ but not both is zero.
\ENDPROOF

\setcounter{theorem}{\value{Lavaurscounter}}
\addtocounter{theorem}{-1}
\THM
Let $0<\eta_-<\eta_+<1$ be rationals for which the parameter rays $R_{\eta_\pm}^\Mbrot$ 
coland at some point $c_0\in\Mbrot$ and let $W_{\eta_-,\eta_+}^\Mbrot$ 
denote the parameter sector bounded by $R_{\eta_-}^\Mbrot\cup\{c_0\}\cup R_{\eta_+}^\Mbrot$ 
and not containing $0$. 
Then the forward orbits of $\eta_\pm$ do not enter the interval $]\eta_-,\eta_+[$. 
And for any $c\in W_{\eta_-,\eta_+}^\Mbrot$ the pair of dynamical rays $R_{\eta_\pm}^c$ move 
homorphically with $c$, co-land at some 
repelling (pre)periodic point $z(c)$ with $Q_c^{k'+l}(z(c)) = Q_c^l(z(c))$, where $l\geq 0$ 
is the common preperiod of $\eta_\pm$ and $k'>0$ divides the common period 
$k>0$ of $\sigma^l(\eta_\pm)$ and the set $R_{\eta_-}^c\cup\{z(c)\}\cup R_{\eta_+}^c$ 
bounds a sector $W^c$ containing $c$, but not $0$.
\ENDTHM

\PROOF\label{proofofLavaurs}
This theorem is at least folklore. We  supply a proof here for completeness.
We shall treat separately the strictly preperiodic case $l>0$ and the periodic case $l=0$. 
For the strictly preperiodic case we have $k=qk'$ with $q>1$ and $c_0$ 
admits precisely $q$ external arguments $0< \theta_0 < \ldots < \theta_{q-1} < 1$ 
both in dynamical plane and in parameter plane by the Douady-Hubbard ray landing theorem. 
The arguments $\eta_-<\eta_+$ are amongst these.
The set 
$$
R^c = \bigcup_{i=0}^{q-1} \overline{R_{\theta_i}^c}
$$ 
moves holomorphically with $c$ in $\C\sm\whR^\Mbrot$, where
$$
\whR^\Mbrot = \bigcup_{i=0}^{q-1}\bigcup_{j=1}^{k+l} \overline{R_{\sigma^j(\theta_i)}^\Mbrot}.
$$
Because the B{\"o}ttcher coordinate $\phi_c$ depends holomorphically on $c$ 
and thus $R^c$ moves holomorphically with $c$ as long as $c$ does not belong to the 
strict forward orbit of $R^c$. 
Write $W_{c_0}^\Mbrot$ for the sector bounded by $\overline{R_{\theta_0}^\Mbrot\cup R_{\theta_{q-1}}^\Mbrot}$. 
Then $W_{\eta_-,\eta_+}^\Mbrot\subseteq W_{c_0}^\Mbrot$ 
and it suffices to prove that $W_{c_0}^\Mbrot\cap\whR^\Mbrot = \emptyset$. 
For the later it is enough to prove that 
\begin{equation}\label{preperiodicdisjointness}
\left(\bigcup_{i=0}^{q-1}\bigcup_{j=1}^{k+l} \sigma^j(\theta_i)\right) \cap [\theta_0,\theta_{q-1}] = \emptyset.
\end{equation}
To this end let us consider the Hubbard tree $T^{c_0}$ for $Q_{c_0}$. 
In this strictly preperiodic case $T^{c_0}$ is the minimal connected subset of $K_c=J_c$ 
containing the orbit 
$$
\OO^{c_0}(0) = \bigcup_{j=0}^{k+l} Q_{c_0}^j(0).
$$
As the orbit $\OO^{c_0}(0)$ is forward invariant, so is $T^{c_0}$. 
Moreover any extremal point of $T^{c_0}$ belongs to $\OO^{c_0}(0)$ by minimality. 
As $Q_ {c_0}^j$ is a local homeomorphism for all $j$ the critical value $c_0 = Q_ {c_0}(0)$ is 
necessarily an extremal point. 
This implies \eqref{preperiodicdisjointness}. 
Notice that the conclusion of the theorem holds in this case even for $c$ in 
a neighbourhood of $\overline{W_{c_0}^\Mbrot}$.

The periodic case is similar and yet slighly different. 
The common landing point $c_0$ of the two parameter rays $R_{\eta_\pm}^\Mbrot$ 
is the root of a hyperbolic component $H\not=H_0$. 
Let us rename $c_0$ to $c_1$ and use $c_0$ to denote the center of $H$. 
As above the dynamical rays $R_{\eta_\pm}^c$ moves holomorphically on 
$\C\sm\whR^\Mbrot$, where
$$
\whR^\Mbrot = \bigcup_{j=0}^{k-1} \overline{R_{\sigma^j(\eta_-)}^\Mbrot\cup R_{\sigma^j(\eta_+)}^\Mbrot}
$$
And to prove the theorem it suffices to prove that 
$R_{\eta_-}^{c_0}$ and $R_{\eta_+}^{c_0}$ coland at a repelling periodic point $z(c_0)$ 
in the dynamical plane of $Q_{c_0}$
and that 
$$ 
\left(\bigcup_{j=0}^{k-1}\sigma^j(\eta_-)\cup\sigma^j(\eta_+)\right) \cap\; ]\eta_-,\eta_+[ \;= \emptyset.
$$ 
Again the proof is that $c_0$ is extremal in the Hubbard tree $T^{c_0}$ for $Q_{c_0}$ 
and that $R_{\eta_\pm}^{c_0}$ coland at a $k'$ periodic point $z(c_0)$, $k' | k$
on the boundary of the Fatou component $F_0$ of $c_0$. 
Notice that in this case the Hubbard tree is defined as the minimal D-H regulated set.  
Where D-H regulated means that for any Fatou component $F$ the image $\phi(F\cap T^{c_0})$ 
under the extended B{\"o}ttcher coordinate consists of radial lines. 
The proof of extremality of $c_0$ in $T_{c_0}$ is the same as in the preperiodic case. 
Also by minimality $z(c_0) = \partial F_0\cap T^{c_0}$ is the unique periodic point on the boundary of $F_0$ 
whose period divides $k$. 
By the Douady-Hubbard ray landing theorem $z(c_0)$ is the common landing point of 
$R_{\eta_\pm}^{c_0}$. 
\ENDPROOF

Addresses:

Carsten Lunde Petersen, IMFUFA,
Roskilde University,
Postbox 260,
DK-4000 Roskilde,
Denmark.
e-mail: lunde@ruc.dk

Pascale Roesch, 
Institut de Math\'ematiques de Toulouse, 
Universit\'e Paul Sabatier,
F-31062 Toulouse Cedex 9
France. 
e-mail: roesch@math.univ-toulouse.fr
\end{document}